\documentclass[preprint]{elsarticle}

\usepackage{lineno,hyperref}
\usepackage{hyperref}
\usepackage{mathtools,amsmath,amssymb,mathabx}

\usepackage{graphicx}
\usepackage{tikz}

\journal{arXiv}

\usepackage{fullpage}

\newtheorem{theorem}{Theorem}
\newproof{pf}{Proof}

\DeclarePairedDelimiter\norm{\lVert}{\rVert}

\makeatletter
\def\th@plain{%
  \thm@notefont{}
  \itshape 
}
\def\th@definition{%
  \thm@notefont{}
  \normalfont 
}
\makeatother
\DeclareFontFamily{U}{matha}{\hyphenchar\font45}
\DeclareFontShape{U}{matha}{m}{n}{ <5> <6> <7> <8> <9> <10> gen * matha <10.95> matha10 <12> <14.4> <17.28> <20.74> <24.88> matha12 }{}
\DeclareSymbolFont{matha}{U}{matha}{m}{n}
\DeclareFontFamily{U}{mathx}{\hyphenchar\font45}
\DeclareFontShape{U}{mathx}{m}{n}{ <5> <6> <7> <8> <9> <10> <10.95> <12> <14.4> <17.28> <20.74> <24.88> mathx10 }{}
\DeclareSymbolFont{mathx}{U}{mathx}{m}{n}
\DeclareMathSymbol{\obot} {2}{matha}{"6B}
\DeclareMathSymbol{\bigobot} {1}{mathx}{"CB}

\newcommand\ra{\rightarrow}
\newcommand\mc{\mathcal}
\newcommand\mb{\mathbb}
\newcommand\mf{\mathbf}
\newcommand\ve{\varepsilon}
\newcommand\xx{\mathfrak{X}}
\newcommand\xf{\mathfrak{F}}
\newcommand\xt{\mathfrak{T}}
\newcommand\rest{\!\restriction}

\DeclareMathOperator{\grad}{grad}
\DeclareMathOperator{\Span}{Span}
\DeclareMathOperator{\Ker}{Ker}
\DeclareMathOperator{\tr}{tr}

\DeclareMathOperator{\id}{Id}

\begin{document}

\begin{frontmatter}

\title{Two-root Riemannian Manifolds}


\author[matf]{Vladica Andreji\'c}
\ead{andrew@matf.bg.ac.rs}

\address[matf]{Faculty of Mathematics, University of Belgrade, Belgrade, Serbia}




\begin{abstract}
Osserman manifolds are a generalization of locally two-point homogeneous spaces.
We introduce $k$-root manifolds in which the reduced Jacobi operator has exactly $k$ eigenvalues. 
We investigate one-root and two-root manifolds as another generalization of locally two-point homogeneous spaces.
We prove that there is no two-root Riemannian manifold of odd dimension.
In twice an odd dimension, we describe all two-root Riemannian algebraic curvature tensors 
and give additional conditions for two-root Riemannian manifolds.
\end{abstract}


\end{frontmatter}


\section{Introduction}

The main feature of the most beautiful and most important Riemannian manifolds
is that they are highly symmetric (they have large groups of isometries).
A connected Riemannian manifold is called a (Riemannian) symmetric space if it has the property that
the geodesic symmetry at any point (it fixes the point and reverses geodesics through that point) 
extends to an isometry of the whole space onto itself.

Symmetric spaces can be observed from plenty different points of view.
For example, they can be locally viewed as Riemannian manifolds for which the curvature tensor is invariant under
all parallel translations.
The algebraic description allowed \'Elie Cartan to develop the theory of symmetric spaces
merged with the theory of semisimple Lie groups which led to a complete classification in 1926 \cite{Ca1,Ca2}.

One refined invariant of a symmetric space is the rank, which is the maximal dimension of a totally geodesic flat submanifold.
The rank is always at least one, with equality when the maximal flat submanifolds are geodesics,
in which case the sectional curvature is positive (compact type) or negative (noncompact type).
Among the Riemannian symmetric spaces, those of rank one are of special importance.

On the other hand, we can consider the cosmological principle which says that the spatial distribution of matter 
in the universe is homogeneous and isotropic at a sufficiently large scale.
A homogeneous Riemannian manifold looks geometrically the same when viewed from any point, 
while an isotropic one has the geometry that does not depend on directions.
A connected Riemannian manifold is called two-point homogeneous if its isometry group is transitive on equidistant pairs of points.
However, a connected Riemannian manifold is isotropic if and only if it is two-point homogeneous, see Wolf \cite[Lemma 8.12.1]{Wo}.

We have a complete classification of these spaces, the compact ones were classified by Wang in 1952 \cite{Wang},
while the noncompact ones by Tits in 1955 \cite{Ti}. 
As a consequence of the classification, it is known that any locally two-point homogeneous Riemannian manifold 
is either flat or locally isometric to a rank one symmetric space, see Helgason \cite[p.535]{He}.

A two-point homogeneous connected Riemannian manifold is isometric to one of the following:
a Euclidean space; a sphere; a real, complex or quaternionic, projective or hyperbolic space;
or the Cayley projective or hyperbolic plane.
More precisely, the classification of these spaces includes:
$\mb{R}^n$, $\mf{S}^n$, $\mb{R}\mf{P}^n$, $\mb{C}\mf{P}^n$, $\mb{H}\mf{P}^n$, $\mb{O}\mf{P}^2$,
$\mb{R}\mf{H}^n$, $\mb{C}\mf{H}^n$, $\mb{H}\mf{H}^n$, and $\mb{O}\mf{H}^2$.
However, note that there are isomorphisms in low dimensions:
$\mb{R}\mf{P}^1\cong \mf{S}^1$, $\mb{C}\mf{P}^1\cong \mf{S}^2$, $\mb{H}\mf{P}^1\cong \mf{S}^4$, $\mb{O}\mf{P}^1\cong \mf{S}^8$,
$\mb{C}\mf{H}^1\cong \mb{R}\mf{H}^2$, $\mb{H}\mf{H}^1\cong \mb{R}\mf{H}^4$, $\mb{O}\mf{H}^1\cong \mb{R}\mf{H}^8$.

Local isometries of a locally two-point homogeneous spaces act transitively 
on the sphere bundle of unit tangent vectors,
and therefore fix the characteristic polynomial of the Jacobi operator there.
In this way we get a generalization of locally two-point homogeneous Riemannian manifolds, called the (globally) Osserman manifolds,
in which the characteristic polynomial (or equivalently, the eigenvalues and their multiplicities)
of a Jacobi operator $\mc{J}_X$ is independent of $X$ from the unit tangent bundle.

The lack of other examples led Osserman \cite{O} to conjecture that the converse might also be true.
The question of whether the converse is true (every Osserman manifold is locally two-point homogeneous)
is known as the Osserman conjecture.
The first results on the Osserman conjecture were given by Chi \cite{Ch}, who
established the affirmative answer for manifolds of dimension not divisible by four.
The largest progress in solving the conjecture was made by Nikolayevsky \cite{Ni4,Ni2,Ni1,Ni3}, 
who proved it in all cases, except the manifolds of dimension $16$ whose reduced Jacobi operator has 
an eigenvalue of multiplicity $7$ or $8$. For recent results regarding the Osserman conjecture one can consult \cite{A9}.

A common way is to consider possible Osserman algebraic curvature tensors that may be realized at a point of a Riemannian manifold.
This algebraic approach brings us a stronger version of the conjecture where instead of globally Osserman manifold 
we have a pointwise Osserman manifold in which the characteristic polynomial of Jacobi operator is independent 
on the unit sphere, but can vary from point to point.

Nikolayevsky followed the approach suggested in \cite{GSV} and showed that the pointwise Osserman condition implies the
existence of a Clifford structure, except in dimension $16$, the only dimension in which there can be Osserman curvature tensors
that do not admit a Clifford structure, as is the case with  $\mb{O}\mf{P}^2$ and $\mb{O}\mf{H}^2$.
In fact, there is no much difference between globally and pointwise Osserman manifolds.
Apart from dimension $16$, where counterexamples of the conjecture could possibly arise, 
the only exceptions are dimensions $2$ and $4$.
In dimension $2$, any Riemannian manifold is pointwise Osserman, while globally Osserman are only those with a constant Gauss curvature.
In dimension $4$, the conjecture for globally Osserman manifolds is proved by Chi \cite{Ch}, but
there exist pointwise Osserman manifolds that are not even locally symmetric, see \cite[Corollary 2.7]{GSV} and \cite{Ol}.

It is worth mentioning the concept of Jacobi-dual Riemannian manifolds that satisfy the so-called Raki\'c duality principle \cite{Ra2},
in which $Y$ is an eigenvector of $\mc{J}_X$ if and only if $X$ is an eigenvector of $\mc{J}_Y$.
It is known that a Riemannian manifold is pointwise Osserman if and only if it is Jacobi-dual \cite{NR1,NR2}.

In this work we consider another generalization.
We say that a Riemannian manifold is $k$-root if the reduced Jacobi operator $\widetilde{\mc{J}}_X$
has exactly $k$ eigenvalues for any nonzero tangent vector $X$.
It is easy to check that two-point homogeneous spaces are one-root or two-root, so they will be the subject of our study.

It is worth noting that the Osserman conjecture for two-root Osserman manifolds was studied in \cite{Ch2}, \cite{GSV}, and \cite[pp.34--35]{GKL}.
In fact, the conjecture for two-root Osserman manifolds depends on the statement announced by Nikolayevsky \cite[Theorem 1.2]{Ni3},
that a two-root Osserman manifold of dimension $16$, with multiplicities $7$ and $8$ is locally isometric to $\mb{O}\mf{P}^2$ or $\mb{O}\mf{H}^2$.
Let us also remark that a connected pointwise two-root Osserman manifold of dimension at least $5$ 
is globally Osserman (see Theorem \ref{globalni}).

It is well known that a connected one-root Riemannian manifold is a space of constant sectional curvature (Theorem \ref{one}).
We prove that in odd dimension there is no two-root Riemannian manifold (Theorem \ref{nepar}).
In this article we describe all two-root Riemannian algebraic curvature tensors in twice an odd dimension (Theorem \ref{prostct}),
and give additional conditions for two-root Riemannian manifolds (Theorem \ref{nepar2}).

\section{Preliminaries}

Let $(M,g)$ be a connected Riemannian manifold of dimension $n\geq 3$.
Its Levi-Civita connection $\nabla$ determines the curvature operator $\mc{R}\in\xt^1_3(M)$ by $\mc{R}(X,Y)=[\nabla_X,\nabla_Y]-\nabla_{[X,Y]}$,
and the associated curvature tensor $R=\mc{R}^{\flat}\in\xt^0_4(M)$ with $R(X,Y,Z,W)=g(\mc{R}(X,Y)Z,W)$ for $X,Y,Z,W\in\xx(M)$.
The underlying geometry of the manifold is related to self-adjoint linear operators $\mc{J}_X\colon \xx(M)\to\xx(M)$ 
called the Jacobi operators that are given by $\mc{J}_X(Y)=\mc{R}(Y,X)X$ for any vector field $X\in\xx(M)$.

It is often convenient to study certain geometric problems in a purely algebraic setting.
Reducing the manifold to an arbitrary point $p\in M$ allows us to deal with 
an algebraic curvature tensor $R_p$ on the scalar product space $(T_pM,g_p)$.

To simplify things, we keep the notation in the following way.
Let $R$ be an algebraic curvature tensor on a (positive definite) scalar product space $(\mc{V},g)$ of dimension $n$,
that is, $R\in\xt^0_4(\mc{V})$ satisfies the usual $\mb{Z}_2$ symmetries as well as the first Bianchi identity.
Let $(E_1,\dots,E_n)$ be an arbitrary orthonormal basis in $(\mc{V},g)$ and let $\ve_X=g(X,X)=\norm{X}^2$ be the squared norm of $X\in\mc{V}$.
Raising the index we have $\mc{R}=R^\sharp\in\xt^1_3(\mc{V})$, so $\mc{R}(X,Y)Z=\sum_{i=1}^n R(X,Y,Z,E_i)E_i$.
As before, the Jacoby operator is a self-adjoint linear operator $\mc{J}_X\colon\mc{V}\to\mc{V}$ 
given by $\mc{J}_X(Y)=\mc{R}(Y,X)X$.

Since $g(\mc{J}_X(Y),X)=0$ and $\mc{J}_X(X)=0$, the Jacobi operator $\mc{J}_X$ for any nonzero $X\in\mc{V}$ 
is completely determined by its restriction $\widetilde{\mc{J}}_X\colon X^\perp \to X^\perp$ called the reduced Jacobi operator.
We are interested in $k$-root curvature tensors in which $\widetilde{\mc{J}}_X$
has exactly $k$ distinct eigenvalues for any nonzero $X\in\mc{V}$.
We say that a Riemannian manifold is $k$-root if its curvature tensor at each point is $k$-root.

The simplest case $k=1$ is associated with one-root manifolds.
If $\widetilde{\mc{J}}_X$ has a single eigenvalue $\ve_X\mu_X$ for $X\neq 0$, then follows $X^\perp=\Ker(\widetilde{\mc{J}}_X-\ve_X\mu_X\id_{X^\perp})$.
For mutually orthogonal nonzero $X,Y\in\mc{V}$ 
we have $\mc{J}_X(Y)=\ve_X\mu_X Y$ and $\mc{J}_Y(X)=\ve_Y\mu_Y X$, so because of $g(\mc{J}_X(Y),Y)=g(\mc{J}_Y(X),X)$ we obtain
a constant sectional curvature $\mu_X=\mu_Y=\kappa(X,Y)$. 
The final touch of the Schur lemma implies the following theorem.

\begin{theorem}\label{one}
A connected one-root Riemannian manifold of dimension $n\geq 3$ is a space of constant sectional curvature.
\end{theorem}

Consider the Jacobi operator $\mc{J}_X$ for a nonzero vector $X=\sum_{i=1}^n x_i E_i$.
The entries of its (symmetric real) matrix related to some orthonormal basis $(E_1,\dots,E_n)$ 
are homogeneous polynomial functions of degree two in coefficients $x_1,\dots,x_n$, 
\begin{equation*}
(\mc{J}_X)_{ab}=g(\mc{J}_X(E_b),E_a)=\sum_{i,j=1}^n R_{bija}x_ix_j.
\end{equation*}
According to the perturbation theory (see Kato \cite[Chapter 2]{Ka}),
the spectrum of $\mc{J}_X$ (unordered $n$-tuple consisting of the repeated eigenvalues) 
depends continuously on $X(x_1,\dots,x_n)$.
However, the $k$-root condition implies no crossing of eigenvalues, 
and consequently the multiplicities of eigenvalues do not change as $X$ varies.
This allows us to label the eigenvalues so that they individually are continuous functions.
Moreover, since the number of distinct eigenvalues of $\mc{J}_X$ is fixed, the eigenvalues of $\mc{J}_X$
depends analytically on the coordinates of $X\neq 0$.

\section{Two-root curvature tensors}

Let $R$ be an algebraic curvature tensor derived at a point of two-root Riemannian manifold.
Then, the reduced Jacobi operator $\widetilde{\mc{J}}_X$ for a nonzero $X\in\mc{V}$ has exactly 
two eigenvalues $\ve_X\mu_X$ and $\ve_X\nu_X$ with constant multiplicities, that is, the characteristic polynomial of Jacobi operator is
\begin{equation*}
\det(\lambda\id-\mc{J}_X)=\lambda(\lambda-\ve_X\mu_X)^p(\lambda-\ve_X\nu_X)^q
\end{equation*}
for fixed integers $p\geq q\geq 1$ with $p+q=n-1$.

Two-point homogeneous Riemannian manifolds as our model spaces are Osserman and have constant both $\mu$ and $\nu$.
It is well known that an Osserman algebraic curvature tensor $R$ is $k$-stein for every $k\in\mb{N}$ (see \cite[Section 1.7]{Gi1}), 
which means that there are constants $C_k$ such that $\tr((\mc{J}_X)^k)=(\ve_X)^kC_k$ holds for any vector $X$.
Thus, a geometric realization of a pointwise Osserman two-root manifold has smooth functions $C_k\in \xf(M)$
with $C_k=p\mu^k+q\nu^k$.
However, we know that $\nabla C_1=0$ for $n\neq 2$ and $\nabla C_2=0$ for $n\not\in\{2,4\}$,
which can be seen in
\cite[pp.164--165]{Be1}, \cite[Theorem 2.4]{GSV}, \cite[Section 1.13, pp.75--78]{Gi1}, and \cite[pp.10--15]{GKL}.
Therefore $C_1$ and $C_2$ are global constants, which allows us to conclude that all $C_k$ are global constants,
and hence we have the following theorem from \cite{A5} (see also \cite[pp.10--16]{GKL}).

\begin{theorem}\label{globalni}
A connected pointwise two-root Osserman manifold of dimension $n\geq 5$ is globally Osserman.
\end{theorem}

Let us focus on the unit sphere $\mc{S}=\{X\in\mc{V} : \ve_X=1\}\subset\mc{V}$.
Consider $\widetilde{\mc{J}}_X-\mu_X\id_{X^\perp}$ as a smooth tangent bundle homomorphism over $\mc{S}$
with the identification $T_X\mc{S}\cong X^\perp$.
Since it has a constant rank $q$, $\Ker(\widetilde{\mc{J}}_X-\mu_X\id_{X^\perp})$ is a smooth subbundle of $T\mc{S}$
(see Lee \cite[Theorem 10.34]{Lee}), that is, a smooth $p$-dimensional distribution on $\mc{S}$.
Similarly, $\Ker(\widetilde{\mc{J}}_X-\nu_X\id_{X^\perp})$ is a smooth $q$-dimensional distribution on $\mc{S}$.

It is well known (see Steenrod \cite[Theorem 27.16]{Stee} and Adams \cite{Ad} for details)
that $\mc{S}$ for $\rho(n)\leq k\leq n-1-\rho(n)$ does not admit a continuous 
$k$-dimensional distribution, where $\rho$ is the Hurwitz--Radon number given by
\begin{equation*}
\rho((2a+1)\cdot 2^{4b+c})=8b+2^c \quad\text{ for }\quad 0\leq c\leq 3,
\end{equation*}
which leaves us with
\begin{equation}\label{adams1}
q<\rho(n).
\end{equation}
The inequality \eqref{adams1} significantly reduces the possibilities for the multiplicities $p$ and $q$.
For example, it immediately removes an odd $n$ because of $\rho(n)=1$, which means that $n$ must be even and proves
the following theorem.

\begin{theorem}\label{nepar}
There is no two-root Riemannian manifold of odd dimension.
\end{theorem}

Without loss of generality we can suppose $\mu_X<\nu_X$, since otherwise we consider $-R$ as a new algebraic curvature tensor.
For any nonzero $X\in\mc{V}$ we define the eigenspaces,
\begin{equation*}
\mc{M}(X)=\Ker(\widetilde{\mc{J}}_X-\ve_X\mu_X\id_{X^\perp}), \quad
\mc{N}(X)=\Ker(\widetilde{\mc{J}}_X-\ve_X\nu_X\id_{X^\perp}),
\end{equation*}
where $\dim\mc{M}(X)=p$ and $\dim\mc{N}(X)=q$, which allows an orthogonal decomposition
\begin{equation*}
X^\perp=\mc{M}(X)\obot \mc{N}(X).
\end{equation*}

For nonzero $X,Y\in\mc{V}$ that satisfy $Y\in\mc{M}(X)$, we can decompose $X=M+N$ such that $M\in\mc{M}(Y)$ and $N\in\mc{N}(Y)$.
Because of
\begin{align*}
g(\mc{J}_X(Y),Y)&=g(\ve_X\mu_XY,Y)=\ve_X\ve_Y\mu_X,\\
g(\mc{J}_Y(X),X)&=g(\ve_Y\mu_YM+\ve_Y\nu_YN,M+N)=\ve_Y\ve_M\mu_Y+\ve_Y\ve_N\nu_Y,
\end{align*}
we have $\ve_X\mu_X=\ve_M\mu_Y+(\ve_X-\ve_M)\nu_Y$, and consequently
\begin{equation}\label{emex}
\ve_M=\ve_X\frac{\nu_Y-\mu_X}{\nu_Y-\mu_Y},
\end{equation}
which gives
\begin{equation*}
0\leq \frac{\nu_Y-\mu_X}{\nu_Y-\mu_Y}\leq 1,
\end{equation*}
and hence $\mu_Y\leq \mu_X\leq \nu_Y$. In a similar fashion, $Y\in\mc{N}(X)$ implies
\begin{equation*}
\ve_N=\ve_X\frac{\mu_Y-\nu_X}{\mu_Y-\nu_Y} \quad \text{and} \quad 0\leq \frac{\nu_X-\mu_Y}{\nu_Y-\mu_Y}\leq 1,
\end{equation*}
and therefore $\mu_Y\leq \nu_X\leq \nu_Y$. Hence

\begin{equation}\label{jedan}
\begin{aligned}
0\neq Y\in\mc{M}(X) \quad &\Longrightarrow \quad \mu_Y\leq \mu_X\leq \nu_Y,\\
0\neq Y\in\mc{N}(X) \quad &\Longrightarrow \quad \mu_Y\leq \nu_X\leq \nu_Y.
\end{aligned}
\end{equation}

The restrictions $\mu\rest_{\mc{S}}\colon \mc{S}\to\mb{R}$ and $\nu\rest_{\mc{S}}\colon \mc{S}\to\mb{R}$ 
are continuous functions on a compact, so their ranges are closed intervals.
Because of $\mc{J}_{tX}/\ve_{tX}=\mc{J}_X/\ve_X$ we obtain $\mu(tX)=\mu(X)$ for all $X\neq 0$ and $t\in\mb{R}$.
Hence, for a nonzero $X\in\mc{V}$ we reach $\mu_X\in[\mu_{\min},\mu_{\max}]$ and $\nu_X\in[\nu_{\min},\nu_{\max}]$, which allows us to define
\begin{equation*}
\mc{U}=\mu^{-1}(\mu_{\min})\cup \{0\}, \quad \mc{W}=\nu^{-1}(\nu_{\max})\cup \{0\}.
\end{equation*}

If $0\neq Y\in\mc{M}(X)$ holds for $0\neq X\in\mc{U}$, then \eqref{jedan} gives $Y\in\mc{U}$,
while \eqref{emex} implies $\ve_M=\ve_X$, that is, $X=M\in\mc{M}(Y)$.
Similarly it can be done for $Y\in\mc{N}(X)$ and $X\in\mc{W}$.
In this way we get some kind of the Raki\'c duality principle when the eigenvalues are extremal,
\begin{equation}\label{dva}
\begin{aligned}
0\neq Y\in\mc{M}(X)\, \wedge\, 0\neq X\in\mc{U} \quad &\Longleftrightarrow \quad 0\neq X\in\mc{M}(Y) \,\wedge\, 0\neq Y\in\mc{U},\\
0\neq Y\in\mc{N}(X)\, \wedge\, 0\neq X\in\mc{W} \quad &\Longleftrightarrow \quad 0\neq X\in\mc{N}(Y) \,\wedge\, 0\neq Y\in\mc{W} .
\end{aligned}
\end{equation}

If we have both $\mc{J}_X(Y)=\ve_X\lambda Y$ and $\mc{J}_Y(X)=\ve_Y\lambda X$ 
for nonzero mutually orthogonal $X,Y\in\mc{V}$ and $\lambda\in\mb{R}$,
then for all $\alpha, \beta\in \mb{R}$, the straightforward calculations (see \cite[Lemma 1]{A4}) gives
\begin{equation}\label{lema1}
\begin{aligned}
\mc{J}_{\alpha X+\beta Y}&(\ve_Y\beta X- \ve_X\alpha Y) =\mc{R}(\ve_Y\beta X- \ve_X\alpha Y, \alpha X+\beta Y) (\alpha X+\beta Y)\\
&=(\ve_X\alpha^2+\ve_Y\beta^2)\mc{R}(X,Y)(\alpha X+\beta Y)\\
&=\ve_{\alpha X+\beta Y}(\beta\mc{J}_Y(X)- \alpha\mc{J}_X(Y))
=\ve_{\alpha X+\beta Y}\lambda (\ve_Y\beta X- \ve_X\alpha Y).
\end{aligned}
\end{equation}

According to \eqref{dva}, $Y\in\mc{M}(X)$ with $X\in\mc{U}$ implies $X\in\mc{M}(Y)$ with $Y\in\mc{U}$,
so \eqref{lema1} yields $\ve_Y\beta X- \ve_X\alpha Y\in\mc{M}(\alpha X+\beta Y)$ with $\alpha X+\beta Y\in\mc{U}$.
Hence, $0\neq U\in\mc{U}$ implies $\Span\{U\}\obot\mc{M}(U)\subseteq\mc{U}$, or consequently $\mc{U}^\perp \subseteq \mc{N}(U)$,
as well as its analogue for $\mc{W}$,
\begin{equation}\label{tri}
\begin{aligned}
0\neq U\in\mc{U} \quad &\Longrightarrow \quad \Span\{U\}\obot\mc{M}(U)\subseteq\mc{U},\\
0\neq W\in\mc{W} \quad &\Longrightarrow \quad \Span\{W\}\obot\mc{N}(W)\subseteq\mc{W}.
\end{aligned}
\end{equation}

Since $\dim(\Span\{U\}\obot\mc{M}(U))=p+1$, $\dim(\Span\{W\}\obot\mc{N}(W))=q+1$, and $\dim\mc{V}=p+q+1>(p+1)+(q+1)$,
the Grassmann formula gives a non-trivial intersection, 
\begin{equation}\label{neprazan}
0\neq(\Span\{U\}\obot\mc{M}(U))\cap(\Span\{W\}\obot\mc{N}(W)) \subseteq \mc{U}\cap\mc{W}.
\end{equation}
The formula \eqref{neprazan} allows to take $0\neq A\in \mc{U}\cap\mc{W}$, as an initial step, 
and exploit its nice properties $\mu_A=\mu_{\min}$ and $\nu_A=\nu_{\max}$.

Due to Theorem \ref{nepar}, $n$ must be even, so we consider the next simplest case of twice an odd dimension $n$.
In that case $\rho(n)=2$, so the inequality \eqref{adams1} gives $q=1$, which means a simple root.
In the following section we consider what happens when one root is simple.

\section{Simple root}

Let us assume that one eigenvalue is simple, that is, $q=1$.
If we suppose $p>n/2$ (which excludes only $n=4$), then according to the Grassmann formula 
any two $\mc{M}$ spaces have a non-trivial intersection.
Thus, for nonzero $X,Y\in\mc{U}$ there exists $0\neq S\in\mc{M}(X)\cap\mc{M}(Y)$, so \eqref{dva} yields $X,Y\in\mc{M}(S)$
with $S\in\mc{U}$, and therefore by \eqref{tri}, $\Span\{X,Y\}\subset\mc{U}$, which proves that $\mc{U}$ is a subspace of $\mc{V}$.

We want to show that $\mu$ is constant, or equivalently $\mc{U}=\mc{V}$.
Assuming the opposite, $\mc{U}\neq\mc{V}$, since $\mc{U}$ is a subspace, applying \eqref{tri} we have $\dim\mc{U}=n-1$, and therefore
\begin{equation}\label{cetiri}
0\neq X\in\mc{U} \quad \Longrightarrow \quad \mc{U}=\Span\{X\}\obot\mc{M}(X) \quad\wedge\quad \mc{U}^\perp=\mc{N}(X).
\end{equation}

Let us start with $0\neq A\in \mc{U}\cap\mc{W}$ from \eqref{neprazan}.
For $0\neq Z\in\mc{N}(A)=\mc{U}^\perp$, by \eqref{dva} we have $A\in\mc{N}(Z)$ with $Z\in\mc{W}$, so $\mc{M}(A)=\mc{M}(Z)=\Span\{A,Z\}^\perp$.
For $0\neq B\in\mc{M}(A)=\mc{M}(Z)$ we have $B\in\mc{U}$, so by \eqref{cetiri}, $Z\in\mc{U}^\perp=\mc{N}(B)$.
Then $g(\mc{J}_B(Z),Z)=g(\mc{J}_Z(B),B)$ gives $\nu_B=\mu_Z=c$.
However, if $B$ and $Z$ are units, then by \eqref{lema1} holds
$\mc{J}_{Z\cos t+B\sin t}(B\cos t-Z\sin t)=c(B\cos t-Z\sin t)$,
which implies $\overline{\mu}(t)=\mu(Z\cos t+B\sin t)=c$ or $\overline{\nu}(t)=\nu(Z\cos t+B\sin t)=c$, for any $t\in\mb{R}$.
The functions $\overline{\mu},\overline{\nu}\colon\mb{R}\to\mb{R}$ are continuous with 
$\overline{\mu}<\overline{\nu}$ and $\overline{\mu}(0)= \overline{\nu}(\pi/2)=c$, so
$\mb{R}=\overline{\mu}^{-1}(c)\sqcup \overline{\nu}^{-1}(c)$ is a disjoint union of non-empty closed sets, which is not possible.

The previous result proves that $\mu_X=\mu$ must be constant for $q=1$, unless $n=4$.
This allows us to introduce a new algebraic curvature tensor $R'=R-\mu R^1$, 
where $R^1\in\xt^0_4(\mc{V})$ is a tensor of constant sectional curvature one given by 
\begin{equation*}
R^1(X,Y,Z,W)=g(Y,Z)g(X,W)-g(X,Z)g(Y,W).
\end{equation*}
This trick shifts the eigenvalues, and the characteristic polynomial of the new Jacobi operator becomes
$\det(\lambda\id-\mc{J}'_X)=\lambda^{n-1}(\lambda-\ve_X(\nu_X-\mu))$.
 
In order not to complicate things too much, we shall keep the previous notation and assume that
$\mc{J}_X$ has a simple eigenvalue $\ve_X\nu_X>0$, while other eigenvalues are all zero.
This essentially means that the original reduced Jacobi operator $\widetilde{\mc{J}}_X$ 
has a simple eigenvalue $\ve_X(\nu_X+\mu)$, while the other root is $\ve_X\mu$ with multiplicity $n-2$.

Let us choose an arbitrary orthonormal basis $(E_1,\dots,E_n)$ in $\mc{V}$.
Then for any nonzero $X=\sum_{i=1}^n x_i E_i\in\mc{V}$, the Jacobi operator $\mc{J}_X$ is of rank one such that
its matrix entries $\mc{J}_{ij}(X)$ are quadratic forms in $n$ variables $x_1,\dots,x_n$.
Any submatrix of order two in a rank one symmetric matrix is singular which gives 
\begin{equation}\label{square1}
\mc{J}_{ii}(X)\mc{J}_{jj}(X)=\mc{J}_{ij}(X)^2
\end{equation}
for all $1\leq i,j\leq n$. 

If we fix some monomial order (for example, the lexicographical order) then there is 
a unique monic (the coefficient of the largest monomial is $1$) $G(X)$ which is the greatest common divisor of all $\mc{J}_{ij}(X)$.
Permuting the basis we can set
\begin{equation*}
\mc{J}_{ii}(X) = \sigma_iG(X)Q_i(X)P_i(X)^2,
\end{equation*}
where $P_i(X)$ and $Q_i(X)$ are some nonzero polynomials for $1\leq i\leq m$, with additional $\mc{J}_{ii}(X)=0$ for $m<i\leq n$,
while $\sigma_i\in\{-1,1\}$.
However, such decomposition is unique up to sign of $P_i(X)$ if we set that $Q_i(X)$ is monic square-free. 
Then
\begin{equation*}
\sigma_i\sigma_jG(X)^2Q_i(X)Q_j(X)P_i(X)^2P_j(X)^2=\mc{J}_{ij}(X)^2 
\end{equation*}
implies $Q_i(X)=Q_j(X)=Q(X)=1$ and $\sigma_i=\sigma_j=\sigma$ for $1\leq i\leq m$, and therefore we have
$\mc{J}_{ij}(X) = \sigma_{ij} G(X)P_i(X)P_j(X)$, where $\sigma_{ij}\in\{-1,1\}$.
Additionally, by \eqref{square1}, $\mc{J}_{ij}(X)=0$ holds whenever $m<i\leq n$ or $m<j\leq n$,
which can be treated as $P_i(X)=0$ for $m<i\leq n$ and extend the indices to $m=n$. 

Another submatrix of order two gives $\mc{J}_{1i}(X)\mc{J}_{ij}(X)=\mc{J}_{1j}(X)\mc{J}_{ii}(X)$, so
$\sigma_{1i}\sigma_{ij}=\sigma_{1j}\sigma_{ii}$.
Because of $\sigma_{ii}=\sigma_i=\sigma$ we have $\sigma_{ij}=\sigma\sigma_{1i}\sigma_{1j}$,
and therefore $\mc{J}_{ij}(X) = \sigma G(X)\sigma_{1i}P_i(X)\sigma_{1j}P_j(X)$.
Since the polynomials $P_i(X)$ are unique up to sign, we can use $\sigma_{1i}P_i(X)$ instead of $P_i(X)$ to obtain
$\mc{J}_{ij}(X) = \sigma G(X)P_i(X)P_j(X)$ for all $1\leq i\leq n$. 

Moreover, comparing the degrees in a polynomial $\mc{J}_{ij}(X)$ we conclude that all $P_i$ have the same degree, zero or one.
The degree zero yields constant polynomials $P_i$, so $\mc{J}_X=G(X)M$, for some constant matrix $M$.
In that case, if $\mc{J}_X(Y)=\ve_X\nu_X Y$, then $\mc{J}_X(Z)=0$ for all $Z\in Y^\perp$, which gives $MZ=0$.
However, then $\mc{J}_Y(Z)=G(Y)MZ=0$, which gives the contradiction $\mc{J}_Y=0$.
Therefore, all $P_i$ have degree one, while $G(X)$ has degree zero and consequently $G(X)=1$. 

Summarizing the previous results, the equation
\begin{equation*}
\mc{J}_{ij}(X) = \sigma P_i(X)P_j(X)
\end{equation*}
holds for all $1\leq i,j\leq n$, where $P_i$ are linear homogeneous polynomials.
If we set 
\begin{equation*}
P(X)=\sum_{i=1}^n P_i(X)E_i,
\end{equation*}
then it follows
\begin{equation*}
\mc{J}_X(P(X)) =\sum_{i=1}^n P_i(X)\sum_{j=1}^n \mc{J}_{ji}(X)E_j =\sigma \sum_{i=1}^n P_i(X)^2 (P(X)). 
\end{equation*}
Thus, $P(X)$ is an eigenvector of $\widetilde{\mc{J}}_X$ associated to the simple eigenvalue 
\begin{equation*}
\sigma \ve_{P(X)}=\sigma \sum_{i=1}^n P_i(X)^2=\tr \mc{J}_X=\ve_X\nu_X,
\end{equation*}
but since we set $\nu_X>0$, it must be $\sigma=1$. In this way we construct a linear map $P\colon\mc{V}\ra\mc{V}$ such that
$\mc{N}(X)=\Span\{P(X)\}$ and 
\begin{equation*}
\nu_X=\frac{\ve_{P(X)}}{\ve_X}
\end{equation*}
with $\nu_{P(X)}\geq \nu_X$ (because of \eqref{jedan}) for any nonzero $X\in\mc{V}$.

Let us start with $0\neq A\in\mc{W}$, when $P(A)\in\mc{N}(A)$, because of \eqref{dva}, implies $A\in\mc{N}(P(A))$ with $P(A)\in\mc{W}$.
Hence, by \eqref{lema1},
\begin{equation*}
\mc{J}_{\alpha A+\beta P(A)} (\ve_{P(A)}\beta A- \ve_A\alpha P(A))
=\ve_{\alpha A+\beta P(A)}\nu_{\max} (\ve_{P(A)}\beta A- \ve_A\alpha P(A)),
\end{equation*}
which gives $\nu(\alpha A+\beta P(A))=\nu_{\max}$ and $P(\alpha A+\beta P(A))\,\propto\, \ve_{P(A)}\beta A- \ve_A\alpha P(A)$.
Using the linearity of $P$ and the fact that $P(A)\perp \Span\{A,P^2(A)\}$, 
we get the coefficient of proportionality equal to $-1/\ve_A$, and consequently $P^2(A)=-\nu_A A$
with $\Span\{A,P(A)\}\subseteq\mc{W}$.

We can continue in a similar manner, using $A_1=A$ and $\nu_1=\nu_{\max}$ as the induction basis.
Let us suppose that we already have mutually orthogonal nonzero vectors $A_1,P(A_1),\dots,A_k,P(A_k)$ such that 
\begin{equation}\label{osobinei}
\Span\{A_i,P(A_i)\}\subseteq\nu^{-1}(\nu_i) \quad \text{and} \quad P^2(A_i)=-\nu_i A_i
\end{equation}
hold for all $1\leq i\leq k$ with $0<\nu_k\leq\dots\leq\nu_1$.
We define  
\begin{equation*}
\nu_{k+1}=\max\{ \nu_X : X\in \mc{S}\cap \mc{M}(A_1)\dots \cap \mc{M}(A_k) \}\leq \nu_k
\end{equation*}
and take arbitrarily $0\neq A_{k+1}\in\nu^{-1}(\nu_{k+1})$.
It is fruitful to notice that, since $\mu$ is constant, the duality \eqref{dva} always provides $Y\in\mc{M}(X)\iff X\in\mc{M}(Y)$.
As a consequence of this, $A_{k+1}\in\mc{M}(A_i)=\mc{M}(P(A_i))=\Span\{A_i,P(A_i)\}^\perp$ implies $A_i,P(A_i)\in\mc{M}(A_{k+1})\perp\mc{N}(A_{k+1})$,
so $P(A_{k+1})\in\Span\{A_i,P(A_i)\}^\perp=\mc{M}(A_i)$.
Thus $\nu_{P(A_{k+1})}\leq \nu_{k+1}$, so by \eqref{jedan} we have $\nu_{P(A_{k+1})}=\nu_{k+1}$, 
while \eqref{lema1} yields $\Span\{A_{k+1},P(A_{k+1})\}\subseteq \nu^{-1}(\nu_{k+1})$ and $P^2(A_{k+1})=-\nu_{k+1}A_{k+1}$.

This procedure uses constants $0<\nu_{n/2}\leq\dots\leq\nu_1$ 
with the properties \eqref{osobinei} to exhaust the space
\begin{equation}\label{razb}
\mc{V}=\bigobot_{i=1}^{n/2} \Span\{A_i,P(A_i)\}.
\end{equation}

Having that on mind, it is easy to see that $P$ is skew-adjoint.
Namely, if we set $X=\sum_{i=1}^{n/2} (x_iA_i+\overline{x}_iP(A_i))$ and $Y=\sum_{i=1}^{n/2} (y_iA_i+\overline{y}_iP(A_i))$, then
\begin{equation*}
g(P(X),Y)=\sum_{i=1}^{n/2} g(x_iP(A_i)-\overline{x}_i\nu_i A_i, y_iA_i+\overline{y}_iP(A_i) )
=\sum_{i=1}^{n/2} \nu_i\ve_{A_i}(x_i\overline{y}_i -\overline{x}_iy_i)=-g(P(Y),X).
\end{equation*}

The key idea is that any skew-adjoint endomorphism $P$ on $\mc{V}$ generates an algebraic curvature tensor $R^P\in\xt^0_4(\mc{V})$ by
\begin{equation*}
(X,Y,Z,W) \mapsto g(PX,Z)g(PY,W) -g(PY,Z)g(PX,W) +2g(PX,Y)g(PZ,W),
\end{equation*}
for all $X,Y,Z,W\in\mc{V}$, which can be easily checked.
Let us remark that these constructions are common for a complex structure $P$ on $(\mc{V},g)$ that preserves the scalar product,
but for our construction the additional condition $P^2=-\id$ is not necessary (see \cite{Gi1} and \cite{AL}).
The corresponding curvature operator has
\begin{equation*}
\mc{R}^P(X,Y)Z= g(PX,Z)PY- g(PY,Z)PX+ 2g(PX,Y)PZ,
\end{equation*}
and consequently the Jacobi operator satisfies 
\begin{equation*}
\mc{J}^P_X(Y)=\mc{R}^P(Y,X)X=3g(PY,X)PX= -3g(Y,PX)PX,
\end{equation*}
that is,
\begin{equation*}
\mc{J}^P_X= \begin{cases}
-3\ve_X\nu_X\id &\text{on } \Span\{P(X)\}\\
0 &\text{on } \Span\{P(X)\}^\perp\end{cases}.
\end{equation*}

Therefore, taking into account the shifting of eigenvalues for $\ve_X\mu$, 
and the possible choice of $\nu_X<\mu_X$ from the beginning of discussion, 
the algebraic curvature tensor must be of form 
\begin{equation}\label{tensor2r}
R=\pm\left(-\frac{1}{3}R^P+\mu R^1\right). 
\end{equation}

This result is better expressed in an orthonormal basis $(E_1,F_1,\dots,E_{n/2},F_{n/2})$ obtained from \eqref{razb} by rescaling
$E_i=A_i/\sqrt{\ve_{A_i}}$ and $F_i=P(A_i)/\sqrt{\ve_{P(A_i)}}$.
Conversely, for any orthonormal basis $(E_1,F_1,\dots,E_{n/2},F_{n/2})$ in $\mc{V}$, constants $0<\nu_{n/2}\leq\dots\leq\nu_1$
define a skew-adjoint endomorphism $P$ on $\mc{V}$ by
\begin{equation}\label{orto}
P(E_i)=\sqrt{\nu_i}F_i,\quad P(F_i)=-\sqrt{\nu_i}E_i,
\end{equation}
for all $1\leq i\leq n/2$.  

\begin{theorem}\label{prostct}
Any two-root algebraic curvature tensor of dimension $n>4$ with a simple root is of the form \eqref{tensor2r},
for $\mu\in\mb{R}$ and some skew-adjoint endomorphism $P$ defined by \eqref{orto} using positive constants $\nu_1,\dots,\nu_{n/2}\in\mb{R}$.
\end{theorem}

\section{Geometric realization}

Theorem \ref{prostct} and the formula \eqref{tensor2r} characterize all possible two-root algebraic curvature tensors of twice an odd dimension.
The second step is then based on the use of the second Bianchi identity with an idea to decide which of these
algebraic curvature tensors may be realized as curvature tensors of a Riemannian manifold.

We shall study the Riemannian manifold $M$ locally in a neighbourhood $U\subset M$ of some point.
There we can set a local orthonormal frame and smoothly extend the elements of our construction.
The smoothness of the curvature tensor $R\in\xt^0_4(U)$ gives the smoothness of $\mu\in\xf(U)$,
while $\nu$ is smooth on the tangent bundle $TU$ minus the zero section.
Then, the way we constructed $P$ brings the smoothness of $P_i(X)\in\xf(U)$, which yields a skew-adjoint operator $P\in\xt^1_1(U)$.
Finally, $\nu\in\xf(TU\setminus(M\times\{0\}))$ implies $\nu_1,\dots,\nu_{n/2}\in\xf(U)$, 
and we can extend our orthonormal bases from the construction to a local orthonormal frame $(E_1,F_1,\dots,E_{n/2},F_{n/2})$ in $\xx(U)$
that fits the formula \eqref{orto}. 
It is convenient to use this frame in the following proof.

Such extensions allow us to apply covariant derivatives to our tensors.
It is important to notice that $\nabla_V P\in\xt^1_1(U)$ is also skew-adjoint,
since $PX\perp X$ implies
\begin{equation*}
0=\nabla_V(g(PX,X))=g(\nabla_V (PX),X)+g(PX,\nabla_VX)=g(\nabla_V (PX),X)-g(P\nabla_VX,X)=g((\nabla_V P)X,X),
\end{equation*}
which after the polarization gives 
\begin{equation*}
g((\nabla_V P)X,Y)=-g((\nabla_V P)Y,X),
\end{equation*}
for all $X,Y,V\in\xx(U)$.

Since $\nabla R^1=0$, the covariant derivative along a vector field $V\in\xx(U)$ 
of our curvature tensor $R$ from the formula \eqref{tensor2r} can be expressed by
\begin{equation*}
\nabla_V R=\mp\frac{1}{3}\nabla_V R^P \pm (\nabla_V\mu) R^1. 
\end{equation*}
For all $X,Y,Z,W,V\in\xx(U)$ we can calculate
\begin{equation}\label{formulica0}
\begin{aligned}
(\nabla_V & R^P)(X,Y,Z,W)\\
=&\phantom{+}g\big( g(PX,Z) (\nabla_VP)Y - g(PY,Z) (\nabla_VP)X  +2g(PX,Y) (\nabla_VP)Z,W\big)\\
&+g\big(g((\nabla_VP)X,Z)PY- g((\nabla_VP)Y,Z)PX +2g((\nabla_VP)X,Y)PZ ,W\big),
\end{aligned}
\end{equation}
and
\begin{equation*}
\begin{aligned}
(\nabla_V & R^P)(X,Y,Y,X)+ (\nabla_X R^P)(Y,V,Y,X)+ (\nabla_Y R^P)(V,X,Y,X)\\
=&\phantom{+} 3g(PX,Y)\big( 2g((\nabla_VP)Y,X) -g((\nabla_XP)Y,V) +g((\nabla_YP)X,V) \big)\\
&-3g\big( g((\nabla_XP)X,Y)PY+ g((\nabla_YP)Y,X)PX , V\big).
\end{aligned}
\end{equation*}
Thus, applying the second Bianchi identity yields
\begin{equation}\label{formulica}
\begin{aligned}
0=&(\nabla_V R)(X,Y,Y,X)+ (\nabla_X R)(Y,V,Y,X)+ (\nabla_Y R)(V,X,Y,X)\\
=&\mp g(PX,Y)\big( 2g((\nabla_VP)Y,X) -g((\nabla_XP)Y,V) +g((\nabla_YP)X,V) \big)\\
&\pm g\big( g((\nabla_XP)X,Y)PY+ g((\nabla_YP)Y,X)PX , V\big) \\
&\pm (\nabla_V \mu) (\ve_X\ve_Y-g(X,Y)^2) \pm (\nabla_X \mu) (g(X,Y)g(Y,V)-\ve_Yg(X,V))\\
&\pm (\nabla_Y \mu) (g(X,Y)g(X,V)-\ve_Xg(Y,V)).
\end{aligned}
\end{equation}
Assuming $Y\perp PX$ in \eqref{formulica} we get
\begin{align*}
0=&\phantom{+} g\big( g((\nabla_XP)X,Y)PY+ g((\nabla_YP)Y,X)PX , V\big) \\
&+(\nabla_X \mu) g(g(X,Y)Y-\ve_YX,V)+(\nabla_Y \mu) g(g(X,Y)X-\ve_XY,V)\\
&+ (\nabla_V \mu) (\ve_X\ve_Y-g(X,Y)^2),
\end{align*}
and therefore
\begin{align*}
(\ve_X\ve_Y & -g(X,Y)^2) (\nabla \mu)^{\sharp}=- g((\nabla_XP)X,Y)PY - g((\nabla_YP)Y,X)PX\\
&+ (X(\mu)\ve_Y-Y(\mu)g(X,Y))X+ (Y(\mu)\ve_X-X(\mu)g(X,Y))Y.
\end{align*}
Thus, for nowhere vanishing $X,Y\in\xx(U)$ such that $Y\in\Span\{X,PX\}^\perp$ we have
$(\nabla \mu)^{\sharp}\in\Span\{X,PX,Y,PY\}$.
However, using our frame with \eqref{orto} we get
\begin{equation*}
(\nabla \mu)^{\sharp}\in\bigcap_{1\leq i<j\leq {n/2}}\Span\{E_i,F_i,E_j,F_j\}=0,
\end{equation*}
which gives $\nabla\mu=0$. Therefore $\mu$ must be constant, while the formula \eqref{formulica} yields 
\begin{equation}\label{formulica2}
\begin{aligned}
& g(PX,Y)\big( 2g((\nabla_VP)Y,X) -g((\nabla_XP)Y,V) +g((\nabla_YP)X,V) \big)\\
&- g\big( g((\nabla_XP)X,Y)PY+ g((\nabla_YP)Y,X)PX , V\big) =0.
\end{aligned}
\end{equation}

Again, $Y\perp PX$ gives $g((\nabla_XP)X,Y)PY+ g((\nabla_YP)Y,X)PX =0$,
while the additional $Y\perp X$ provides linear independence for $PX$ and $PY$ (since $X$ and $Y$ are linearly independent as orthogonal),
and therefore $g((\nabla_XP)X,Y)=0$ for $Y\in\Span\{X,PX\}^\perp$.
However, we already know that $g((\nabla_XP)X,X)=0$, which implies $(\nabla_XP)X\,\propto\, PX$.

Let us define the map $\lambda_X\colon U\to \mb{R}$ for any $X\in\xx(U)$ by
$\lambda_X=g((\nabla_XP)X, PX)/g(PX,PX)$ on $(\ve_{PX})^{-1}(\mb{R}_+)=(\ve_{X})^{-1}(\mb{R}_+)\subseteq U$
and $\lambda_X=0$ on $(\ve_X)^{-1}(\{0\})$, where the previously proven proportionality yields
\begin{equation}\label{prop}
(\nabla_XP)X=\lambda_X PX.
\end{equation}
It is easy to check that $\lambda_{f X}=f\lambda_X$ holds for $f\in\xf(U)$.
On the other hand, from \eqref{prop} for any $X,Y\in\xx(U)$ we have
\begin{equation*}
(\nabla_X P)Y + (\nabla_Y P)X=(\lambda_{X+Y}-\lambda_X)PX+(\lambda_{X+Y}-\lambda_Y)PY,
\end{equation*}
which after taking the inner product by $X$ gives
\begin{equation*}
(\lambda_{X+Y}-\lambda_Y) g(PY,X)=g((\nabla_X P)Y,X)=-g((\nabla_X P)X,Y)=-\lambda_X g(PX,Y),
\end{equation*}
and therefore $(\lambda_{X+Y}-\lambda_Y-\lambda_X) g(PX,Y)= 0$. Hence, the additivity
\begin{equation*}
\lambda_{X+Y}=\lambda_X+\lambda_Y
\end{equation*}
holds whenever $g(Y,PX)$ is nowhere zero. 
For any $p\in U$, the condition $Y_p\perp PX_p$ can be excluded by continuity of $\lambda(p)\colon T_pU\to \mb{R}$
given by $\lambda(p)(X_p)=\lambda_X(p)$.
Thus, the additivity holds for all $X,Y\in\xx(U)$, which means that $\lambda$ is $\xf(U)$-linear.
Consequently, since $\lambda_E\in\xf(U)$ for a unit $E\in\xx(U)$, we have $\lambda_X\in\xf(U)$ for any $X\in\xx(U)$,
and finally $\lambda\in\xt^0_1(U)=\xx^*(U)$.

With this in mind, the equation \eqref{formulica2} becomes
\begin{equation*}
g(PX,Y) \Big( 2g((\nabla_VP)Y,X) -g\big((\nabla_XP)Y - (\nabla_YP)X + \lambda_X PY - \lambda_Y PX ,V\big) \Big)=0.
\end{equation*}
Hence,
\begin{equation}\label{formulica3}
2g((\nabla_VP)Y,X) = g\big((\nabla_XP)Y - (\nabla_YP)X + \lambda_X PY - \lambda_Y PX ,V\big)
\end{equation}
holds in the case that $g(PX,Y)$ is nowhere zero. However, since the right hand side is linear in $Y$
and there is a frame consisting of vector fields that are not orthogonal to $PX$, 
the equation \eqref{formulica3} holds for all $X,Y,V\in\xx(U)$.
Applying \eqref{formulica3} twice, we have
\begin{equation*}
\begin{aligned}
4g((\nabla_VP)Y,X) &= g\big((\nabla_VP)Y - (\nabla_YP)V + \lambda_V PY - \lambda_Y PV ,X\big)\\
&+ 2g\big(- (\nabla_YP)X + \lambda_X PY - \lambda_Y PX ,V\big),
\end{aligned}
\end{equation*}
and therefore
$2\lambda_X g(PY,V)= g\big(3(\nabla_VP)Y -(\nabla_YP)V - \lambda_V PY - \lambda_Y PV , X \big)$,
which gives
\begin{equation*}
2g(PY,X)\lambda^\sharp  =  3(\nabla_XP)Y -(\nabla_YP)X - \lambda_X PY - \lambda_Y PX.
\end{equation*}

On the other hand, from the definition of $\lambda\in\xx^*(U)$ we have $(\nabla_XP)Y+(\nabla_YP)X=\lambda_Y PX+\lambda_X PY$, and therefore
we obtain $g(X,PY)\lambda^\sharp=(\nabla_XP)Y-(\nabla_YP)X$, which can be written as
\begin{equation}\label{lambda}
2(\nabla_XP)Y =g(X,PY) \lambda^\sharp  +\lambda_Y PX + \lambda_X PY.
\end{equation}

Now that we know $\nabla P$, it remains to calculate $\nabla^2 P$ and use the Ricci identity for the tensor field $P\in\xt^1_1(U)$.
From \eqref{lambda} we have
\begin{equation*}
\begin{aligned}
2(\nabla_X\nabla_Y P)Z=& 2\nabla_X ((\nabla_Y P)Z)- 2(\nabla_Y P)(\nabla_X Z)\\
=&\nabla_X ( g(Y,PZ) \lambda^\sharp  +\lambda_Z PY + \lambda_Y PZ) 
- (g(Y,P(\nabla_X Z)) \lambda^\sharp  +\lambda_{\nabla_X Z} PY + \lambda_Y P(\nabla_X Z))\\
=& (g(\nabla_XY,PZ)+g(Y,\nabla_X PZ)) \lambda^\sharp + g(Y,PZ)\nabla_X \lambda^\sharp\\
&+ ( g(\nabla_X\lambda^\sharp,Z) + g(\lambda^\sharp,\nabla_XZ) )PY +\lambda_Z \nabla_X PY \\
&+ ( g(\nabla_X\lambda^\sharp,Y) + g(\lambda^\sharp,\nabla_XY) )PZ +\lambda_Y \nabla_X PZ \\
&- g(Y,P(\nabla_X Z)) \lambda^\sharp  -\lambda_{\nabla_X Z} PY - \lambda_Y P(\nabla_X Z)\\
=& (g(\nabla_XY,PZ)+g(Y,(\nabla_X P)Z)) \lambda^\sharp +g(\nabla_X\lambda^\sharp,Z) PY+ ( g(\nabla_X\lambda^\sharp,Y) + g(\lambda^\sharp,\nabla_XY) )PZ\\
& +g(Y,PZ)\nabla_X \lambda^\sharp +\lambda_Z \nabla_X PY +\lambda_Y (\nabla_X P)Z, 
\end{aligned}
\end{equation*}
and therefore
\begin{equation*}
\begin{aligned}
2(\nabla^2_{X,Y}P&-\nabla^2_{Y,X}P)Z=2(\nabla_X\nabla_Y P-\nabla_Y \nabla_X P-\nabla_{\nabla_X Y-\nabla_Y X}P)Z\\
=& (g(\nabla_XY,PZ)+g(Y,(\nabla_X P)Z)) \lambda^\sharp +g(\nabla_X\lambda^\sharp,Z) PY+ ( g(\nabla_X\lambda^\sharp,Y) + g(\lambda^\sharp,\nabla_XY) )PZ\\
& +g(Y,PZ)\nabla_X \lambda^\sharp +\lambda_Z \nabla_X PY +\lambda_Y (\nabla_X P)Z \\
&-(g(\nabla_YX,PZ)+g(X,(\nabla_Y P)Z)) \lambda^\sharp -g(\nabla_Y\lambda^\sharp,Z) PX -( g(\nabla_Y\lambda^\sharp,X) + g(\lambda^\sharp,\nabla_YX) )PZ\\
& -g(X,PZ)\nabla_Y \lambda^\sharp -\lambda_Z \nabla_Y PX -\lambda_X (\nabla_Y P)Z \\
& -g({\nabla_X Y-\nabla_Y X},PZ) \lambda^\sharp  -\lambda_Z P(\nabla_X Y-\nabla_Y X) - \lambda_{\nabla_X Y-\nabla_Y X} PZ,\\
=& (g(Y,(\nabla_X P)Z)-g(X,(\nabla_Y P)Z)) \lambda^\sharp -g(\nabla_Y\lambda^\sharp,Z) PX +g(\nabla_X\lambda^\sharp,Z) PY\\ 
&+ (g(\nabla_X\lambda^\sharp,Y)- g(\nabla_Y\lambda^\sharp,X))PZ +g(Y,PZ)\nabla_X \lambda^\sharp -g(X,PZ)\nabla_Y \lambda^\sharp \\
&+\lambda_Z (\nabla_X P)Y -\lambda_Z (\nabla_Y P)X +\lambda_Y (\nabla_X P)Z-\lambda_X (\nabla_Y P)Z. 
\end{aligned}
\end{equation*}
Applying \eqref{lambda} again we obtain
\begin{equation*}
\begin{aligned}
2(\nabla^2_{X,Y}P&-\nabla^2_{Y,X}P)Z\\
=& \lambda_Z g(Y,PX)\lambda^\sharp -g(\nabla_Y\lambda^\sharp,Z) PX +g(\nabla_X\lambda^\sharp,Z) PY\\
&+(g(\nabla_X\lambda^\sharp,Y)- g(\nabla_Y\lambda^\sharp,X))PZ +g(Y,PZ)\nabla_X \lambda^\sharp -g(X,PZ)\nabla_Y \lambda^\sharp \\
&+\lambda_Z g(X,PY) \lambda^\sharp +\frac{1}{2}(\lambda_Y g(X,PZ)-\lambda_X g(Y,PZ)) \lambda^\sharp 
+\frac{1}{2}\lambda_Y\lambda_Z PX - \frac{1}{2}\lambda_X\lambda_Z PY\\ 
=&\frac{1}{2}(\lambda_Y g(X,PZ)-\lambda_X g(Y,PZ)) \lambda^\sharp 
+(\frac{1}{2}\lambda_Y\lambda_Z -g(\nabla_Y\lambda^\sharp,Z)) PX
-(\frac{1}{2}\lambda_X\lambda_Z -g(\nabla_X\lambda^\sharp,Z)) PY\\
&+(g(\nabla_X\lambda^\sharp,Y)- g(\nabla_Y\lambda^\sharp,X))PZ
+g(Y,PZ)\nabla_X \lambda^\sharp -g(X,PZ)\nabla_Y \lambda^\sharp. 
\end{aligned}
\end{equation*}

To simplify the notation we introduce the operator $Q\in\xt^1_1(U)$ defined by $QX=\frac{1}{2}\lambda_X\lambda^\sharp-\nabla_X \lambda^\sharp$,
so the previous equation becomes
\begin{equation}\label{2rootq}
\begin{aligned}
2(\nabla^2_{X,Y}P-\nabla^2_{Y,X}P)Z=&g(X,PZ)QY-g(Y,PZ)QX \\
&+g(Z,QY)PX -g(Z,QX)PY+ (g(QY,X)-g(QX,Y))PZ.
\end{aligned}
\end{equation}

On the other hand, for the curvature operator $\mc{R}$ from \eqref{tensor2r} we have $\pm 3\mc{R}=-\mc{R}^P+3\mu \mc{R}^1$, 
and therefore
\begin{equation*}
\begin{aligned}
\pm 3(\mc{R}(X,Y)PZ &-P(\mc{R}(X,Y)Z))\\
=& -\mc{R}^P (X,Y)PZ+3\mu \mc{R}^1(X,Y)PZ +P(\mc{R}^P(X,Y)Z) -3\mu P(\mc{R}^1(X,Y)Z)\\
=& -g(PX,PZ)PY +g(PY,PZ)PX -2g(PX,Y)P^2Z +3\mu (g(Y,PZ)X-g(X,PZ)Y)\\
&+g(PX,Z)P^2Y- g(PY,Z)P^2X+ 2g(PX,Y)P^2Z -3\mu (g(Y,Z)PX-g(X,Z)PY)\\
=& -g(PX,PZ)PY +g(PY,PZ)PX  +3\mu g(Y,PZ)X -3\mu g(X,PZ)Y\\
&+g(PX,Z)P^2Y- g(PY,Z)P^2X -3\mu g(Y,Z)PX +3\mu g(X,Z)PY.
\end{aligned}
\end{equation*}
We introduce the self-adjoint operator $S\in\xt^1_1(U)$ defined by $SX=3\mu X+P^2X$, so the previous equation becomes
\begin{equation*}
\pm 3(\mc{R}(X,Y)PZ-P(\mc{R}(X,Y)Z))= g(PX,Z)SY -g(PY,Z)SX +g(SX,Z)PY-g(SY,Z)PX.
\end{equation*}
The Ricci identity 
\begin{equation*}
((\nabla_X\nabla_Y -\nabla_Y \nabla_X -\nabla_{[X,Y]})P)Z= (\nabla^2_{X,Y}P-\nabla^2_{Y,X}P)Z=\mc{R}(X,Y)PZ-P(\mc{R}(X,Y)Z) 
\end{equation*}
holds for all $X,Y,Z\in\xx(U)$, which using \eqref{2rootq} yields
\begin{equation*}
\begin{aligned}
& g(X,PZ)QY-g(Y,PZ)QX +g(Z,QY)PX -g(Z,QX)PY+ (g(QY,X)-g(QX,Y))PZ\\
=&\pm\frac{2}{3}\Big( g(PX,Z)SY -g(PY,Z)SX +g(SX,Z)PY-g(SY,Z)PX\Big).
\end{aligned}
\end{equation*}

It is convenient to introduce another operator $K=Q\pm\frac{2}{3}S \in\xt^1_1(U)$, that is 
\begin{equation*}
KX=\frac{1}{2}\lambda_X\lambda^\sharp-\nabla_X \lambda^\sharp \pm\frac{2}{3}(3\mu X+P^2X),
\end{equation*}
for $X\in\xx(U)$, so the previous equation becomes 
\begin{equation*}
g(X,PZ)KY -g(Y,PZ)KX +g(Z,KY)PX-g(Z,KX)PY +(g(KY,X)-g(KX,Y))PZ=0.
\end{equation*}
The special case $Z=Y$ implies 
\begin{equation}\label{riccik}
g(X,PY)KY  +g(Y,KY)PX +(g(KY,X)-2g(KX,Y))PY=0,
\end{equation}
which holds for all $X,Y\in\xx(U)$.
For an arbitrary nowhere vanishing $Y\in\xx(U)$ we can take a nowhere vanishing $X\in\Span\{Y,PY\}^{\perp}$. In this case $X\perp PY$ gives 
$g(Y,KY)PX +(g(KY,X)-2g(KX,Y))PY=0$, but since $X$ and $Y$ are linearly independent as mutually orthogonal, 
$PX$ and $PY$ are linearly independent, which implies $g(Y,KY)=0$.

Hence, $KY\perp Y$ holds for any $Y\in\xx(U)$, which after the polarization gives $g(KX,Y)+g(KY,X)=0$, and proves that $K$ is also skew-adjoint.
With this in mind, the equation \eqref{riccik} becomes
\begin{equation*}
g(X,PY)KY  +3g(KY,X) PY=0,
\end{equation*}
and holds for all $X,Y\in\xx(U)$.
Substituting $X=PY$ for a nowhere vanishing $Y\in\xx(U)$, we obtain $KY\propto PY$,
while taking the inner product by $PY$ we get $4\ve_{PY}g(KY,PY)=0$, and therefore $K=0$. 
Thus arises the important formula 
\begin{equation}\label{ricci2}
\nabla_X \lambda^\sharp=\frac{1}{2}\lambda_X\lambda^\sharp \pm\frac{2}{3}(3\mu X+P^2X).
\end{equation}

For any $X,Z\in\xx(U)$, we use \eqref{lambda} to calculate
\begin{equation*}
\begin{aligned}
\nabla_X(\ve_{PZ})&=\nabla_Xg(PZ,PZ)=2g(\nabla_X(PZ),PZ)=2g((\nabla_X P)Z+P\nabla_X Z,PZ)\\
&=g(X,PZ) \lambda_{PZ}  + g(PX,PZ) \lambda_Z + g(PZ,PZ)\lambda_X-2g(\nabla_X Z,P^2Z). 
\end{aligned}
\end{equation*}
On the other hand, $\nabla_X(\ve_{PZ})= \nabla_X(\nu_Z\ve_Z)=\nu_Z\nabla_X\ve_Z+\ve_Z\nabla_X\nu_Z$, which gives
\begin{equation*}
\nu_Z\nabla_X\ve_Z+\ve_Z\nabla_X\nu_Z=g(X,PZ) \lambda_{PZ}  - g(X,P^2Z) \lambda_Z + \ve_{PZ}\lambda_X-2g(\nabla_X Z,P^2Z). 
\end{equation*}
Consider the eigenspaces of $P^2$, defined by 
$\mc{P}_j=\ker(P^2+\nu_j\id)=\bigobot_{\nu_i=\nu_j} \Span\{E_i,F_i\}$.
If we suppose that $Z\in\mc{P}_j$ holds for some $1\leq j\leq n/2$, then 
$-2g(\nabla_X Z,P^2Z)=2\nu_j g(\nabla_X Z,Z)=\nu_j \nabla_X \ve_Z$, which implies
\begin{equation*}
\ve_Z\nabla_X\nu_j = g(X,PZ)\lambda_{PZ} +\nu_j g(X,Z) \lambda_Z + \nu_j\ve_{Z}\lambda_X. 
\end{equation*}
Hence, we obtain
\begin{equation}\label{lambda1}
d(\ln \nu_j)(X)=\frac{\nabla_X\nu_j}{\nu_j} 
=\begin{cases}2\lambda_X \text{ for } X\in\Span\{Z,PZ\}\\ \lambda_X \text{ for } X\in \Span\{Z,PZ\}^{\perp}\end{cases}, 
\end{equation}
whenever $Z\in\mc{P}_j$, and therefore
\begin{equation}\label{lambda2}
\lambda=\frac{2}{n+2}\sum_{i=1}^{n/2} d(\ln \nu_i)=\frac{2}{n+2}d(\ln(\nu_1\nu_2\cdots\nu_{n/2})).
\end{equation}

The equation \eqref{lambda2} shows that $\lambda$ cannot be any covector field, but at least one that is the differential of a smooth function.
Moreover, using \eqref{ricci2} we have the necessary condition,
\begin{equation*}
\nabla_X \grad(\ln(\nu_1\nu_2\cdots\nu_{n/2}))
=\frac{1}{2}\lambda_X \grad(\ln(\nu_1\nu_2\cdots\nu_{n/2})) \pm\frac{n+2}{3}(3\mu X+P^2X),
\end{equation*}
that holds for any $X\in\xx(U)$.

\begin{theorem}\label{nepar2}
A two-root Riemannian manifold of dimension $n\equiv 2\pmod 4$ locally has the curvature tensor of the form \eqref{tensor2r},
for a constant $\mu$ and some skew-adjoint linear operator $P$ defined by \eqref{orto} 
using positive smooth functions $\nu_1,\dots,\nu_{n/2}$.
In addition, the equations \eqref{lambda}, \eqref{ricci2}, \eqref{lambda1}, and \eqref{lambda2} hold. 
\end{theorem}

The most natural case has $\lambda=0$, where the equation \eqref{lambda} implies $\nabla P=0$, so \eqref{formulica0} 
gives $\nabla R^P=0$, and consequently $\nabla R=0$, which means that $M$ is locally symmetric.
Moreover, the equation \eqref{ricci2} for $\lambda=0$ implies $P^2=-3\mu\id$, which implies $\nu=3\mu$,
and consequently $M$ is globally Osserman, where the reduced Jacobi operator $\widetilde{\mc{J}}_X$ has
a simple eigenvalue $4\ve_X\mu$, while the other eigenvalue (with multiplicity $n-2$) is four times smaller.
Thus, a connected two-root Riemannian manifold of dimension $n\geq 3$ with $n\equiv 2\pmod 4$ that has $\lambda=0$ 
is globally Osserman, and hence is two-points homogeneous.

Let us remark, that if $\nabla \nu_j=0$ holds for some $1\leq j\leq n/2$, then $\lambda=0$, and the previous conclusion holds.
The question whether there are two-root Riemannian manifolds of twice an odd dimension that are not Osserman remains open
and requires a construction of concrete manifolds with $\lambda\neq 0$.
Let us remark that the first attempt could be $\lambda^\sharp=E_k$ for some $1\leq k\leq n/2$, where \eqref{lambda1} yields
$d(\ln \nu_k)=2\lambda=2d(\ln \nu_i)$ for any $i\neq k$, and therefore there exist constants $C_i$ such that $\nu_k=C_i\nu_i^2$.

\section*{Acknowledgements}

The author was partially supported by the Ministry of Education, Science and Technological
Developments of the Republic of Serbia: grant number 451-03-68/2022-14/200104.

The author is extremely thankful to the referee for valuable suggestions,
especially for the arguments in the final section that shortened some parts of the paper.

\end{document}